

\def\d        #1#2{#1_#2}
\def\dd       #1#2#3{#1_{#2#3}}

\def\dddd     #1#2#3#4#5{#1_{#2#3#4#5}}
\def\u        #1#2{#1^#2}
\def\uu       #1#2#3{#1^{#2#3}}

\def\ud       #1#2#3{{#1}^{#2}_{~#3}}
\def\du       #1#2#3{#1_#2^{~#3}}

\def\dudd     #1#2#3#4#5{#1_{#2~#4#5}^{~#3}}

\def\td       #1#2#3{{d {#1}\over d #2_{#3}}}

\def\cd       #1#2{\nabla_#2#1}

\def\cdb      #1#2{\overline{\nabla}_#2#1}

\def\t        #1{\widetilde {#1}}

\def\b        #1{\overline{#1}}
\def\real     #1{{\rm I\hskip -2pt R^{#1}}}

\def\one      #1{{#1}\acute{}~}
\def\two      #1{{#1}\acute{}\hskip 2pt\acute{}~}

\def\onesub   #1#2{{#1}_{#2}\hskip -2pt\acute{}~}
\def\twosub   #1#2{{#1}_{#2}\hskip -2pt\acute{}\hskip 2pt\acute{}~}

\def\ex       #1#2{\exp_{\vert {#1}}\hskip -2pt {#2}}

\def\lemma    #1{\bigskip\noindent\it{$\underline{\bf Lemma~#1:}$}\medskip}
\def\theorem  #1{\bigskip\noindent\it{$\underline{\bf Theorem~#1:}$}\medskip}
\def\proof    {\bigskip\noindent\item{Proof:~}\rm}
\def\corollar #1{\bigskip\noindent\it{$\underline{\bf Corollary~#1:}$}\medskip}
\def\formula  #1#2{$$#2\leqno#1$$}
\def\remark   {\bigskip\noindent{\bf Remark:~}}
\def\abstract #1{{\narrower\narrower\noindent\klein ABSTRACT. #1
              \par\par}}
\def\section  #1{\vskip 1truecm\noindent{\bf #1}\bigskip\noindent}
\def\prop     #1{\bigskip\noindent\it{$\underline{\bf Proposition~#1:}$}
              \medskip}
\def\cite     #1#2#3#4#5#6{\vskip 0pt\noindent\item {[#1]}{\rm #2}{\it ~#3}
              {\klein #4~#5~#6}\vskip 0pt\noindent}
\def\defin    #1{\bigskip\noindent\it{$\underline{\bf
Definition~#1:}$}\medskip}

\font\klein=cmr8


\def\date{May 02, 1996}
\magnification=\magstep1
\baselineskip =13pt
\noindent
\centerline{\bf    THE FERMI FLOW AND ITS APPLICATION TO GEOMETRY}
\bigskip\noindent
\centerline{\rm Knut Smoczyk$~^{1~2}$}
\bigskip\noindent
\centerline{\rm Mathematics Department, Harvard University, Cambridge,
MA 02138, USA}
\bigskip\noindent
\centerline{\rm \date}
\footnote{}{$^1$ Supported by a Feodor Lynen Fellowship of the Alexander
von Humboldt Foundation}
\footnote{}{$^2$ email: ksmoczyk@abel.math.harvard.edu}
\bigskip\noindent
\abstract{We introduce the notion of Fermi flow for hypersurfaces in
Riemannian manifolds. It turns out that this is a powerful tool to
study the geometry of distance surfaces about a given initial hypersurface.
Some of the results in this paper are known in one form or another, however
the aim is to demonstrate how they can all be derived by the
same method and proved in a very simple manner. In addition we obtain some
new results and results that are stronger than those stated in the literature.}

\section {1. Introduction and preliminary results}
\smallskip\noindent
Let $~N^{n+1}~$ be a Riemannian manifold and $~M^n~$ be an orientable
manifold (possibly with boundary) that is smoothly immersed
by a  local diffeomorphism $~F~$ into $~N^{n+1}~$. Further assume that $~\nu~$
is a normal vector field on $~M_0:=F(M^n)~$. Then we say that $~M^n~$
evolves under the Fermi flow if there exists a
smooth family of local diffeomorphisms $~F_t:M^n\to M_t:=F_t(M^n)\subset
N^{n+1}~$
such that
\formula{(1)}{{d\over dt}F_t(x)=\nu(x,t)~,~F_0=F~.}
If $~M^n~$ is compact, then a smooth solution of this nonlinear first order
PDE exists at least on some short time interval. Indeed the solutions of (1)
are simply given by the foliation arising from Fermi coordinates about
$~M_0~$ since the surfaces are moving by constant speed in their normal
direction (see section 2 for details).
\remark
A more familiar expression for Fermi coordinates in codimension $1$
is Gaussian normal coordinates. But the notation ``Gaussian flow" would
be rather confusing, since this could be mistaken for the flow defined
by the Gaussian curvature on $~M~$.

\smallskip\noindent
If $~\t F:\tilde M^n\to N^{n+1}~$ is another immersion such that
$~\t F(y)=F(x),~\tilde\nu(y)=\nu(x)~$ for a fixed
pair $~(y,x)\in \tilde M^n\times M^n~$
then $~\tilde F_t(y)=F_t(x)~$ as long as a smooth solution of (1) exists.
Roughly speaking, this means that the solutions of (1)
depend only on the exterior geometry.
In a more general setting, in [S], evolution equations have been derived
for various curvature flows given by speed
functions $~-f~$.
Here we have $~f=-1~$. Usually one allows $~f~$ to depend on the second
fundamental form of $~M_t~$. Here it does not, so ``curvature" function is
somewhat misleading. However, this flow measures the effect of exterior
curvature on the motion of hypersurfaces.
In the forthcoming we will assume that $~M^n~$ is compact and we also
set $~r=t~$, since
$$M_t=\{y\in N\vert y{\rm ~can~be~joined~with~}
M_0{\rm~by~a~geodesic~of~length~}t{\rm ~with~initial~speed~}\nu\}$$
as long as a smooth solution of (1) exists.
\medskip\noindent
Let $~\bar g_{\alpha\beta}~$ be
the metric on $~N~$, $~\sigma_{ij}~$, $~\tau_{ij}~$
be the pullbacks of the metric
and second fundamental form on $~M_r~$ (e.g. $~\sigma_{ij}=\bar g_{\alpha\beta}
{\partial F^\alpha\over\partial x^i}{\partial F^\beta\over\partial x^j}~$).
Latin indices range from $~1~$ to $~n~$, greek from $~0~$ to $~n~$ and we
will always sum over repeated indices. In particular to avoid confusion
let us denote the initial pullbacks of the induced metric and second
fundamental
form on $~M_0~$
with $~\rho_{ij},~\lambda_{ij}~$ resp. Upper indices are always meant (if not
otherwise mentioned) with respect to the metric $~\sigma_{ij}~$, e.g.
$~\ud {\tau}ik=\uu {\sigma}ij\dd {\tau}jk~$ where as usual upper indices
for metrics denote the inverse.
\medskip\noindent
For $~1\le k~$ we define
$$A_k:=\ud {\tau}{i_1}{i_2}\cdots \ud {\tau}{i_{k-2}}{i_{k-1}}
\ud {\tau}{i_{k-1}}{i_k}\ud {\tau}{i_k}{i_1}$$
Then the calculations in [S] immediatly imply the following evolution
equations

{\lemma {1.1}\itemitem{(a)} $\td {}r{} \dd {\sigma}ij=2\dd {\tau}ij$
\smallskip\itemitem {(b)} $\td {}r{} d\mu=\tau d\mu$
\smallskip\itemitem {(c)} $\td {}r{}\nu=0$
\smallskip\itemitem {(d)} $\td {}r{} \dd {\tau}ij=\du {\tau}ik\dd {\tau}kj
-\dddd {\b R}0i0j$
\smallskip\itemitem {(e)} $\td {}r{}\ud {\tau}ij=-\ud {\tau}il\ud {\tau}lj
-\dudd {\b R}0i0j$
\smallskip\itemitem {(f)} $\td {}r{}A_k=-kA_{k+1}-k\ud {\tau}{i_1}{i_2}
\cdots \ud {\tau}{i_{k-1}}{i_k}\dudd {\b R}0{i_1}0{i_k}$
}
\bigskip\noindent
$\tau~$ is the mean curvature, $~d\mu~$ the volume form
on $~M_r~$ and $~\b R~$ is the Riemann curvature tensor on $~N~$.

\bigskip\noindent
In the following let $~\d r0~$ denote the maximal distance (time) for which a
smooth solution of (1) exists. We will abbreviate derivatives with respect to
$~r~$ with a prime.
Now for any fixed $~x~$ let $~\eta_{ij}(x,r)~$
be the solution of the pointwise second order ODE
$$\twosub {\eta}{ij}+\dddd {\b R}0j0k\eta^{kl}\rho_{li}=0$$
$$\eta_{ij}(0)=\rho_{ij}~~~;~~~\onesub {\eta}{ij}(0)=\lambda_{ij}~.$$
Here $~\eta^{kl}~$ denotes the inverse of $~\eta_{ij}~$.
Note that there exists a unique, invertible solution up to some distance since
$~\eta_{ij}(0)~$ is invertible.
Then we have

{\theorem{1.2}\itemitem{(a)} $\sigma
_{ij}=\rho^{kl}\eta_{ki}\eta_{lj}~~;~~\sigma^{mn}=\eta^
{ms}\eta^{nt}\rho_{st}$
\smallskip\itemitem{(b)}  $\tau_{ij}=\rho^{kl}\onesub {\eta}
{ki}\eta_{lj}=\rho^{kl}\eta_{ki}\onesub {\eta}{lj}$
\smallskip\itemitem{(c)} $d\mu_t=\det(\rho^{ik}\eta_{kj})d\mu_0~.$
}
\proof
$\sigma_{ij}~$ and $~a_{ij}:=\rho^{kl}\eta_{ki}\eta_{lj}~$ both satisfy the
pointwise second order ODE
$$\twosub y{ij}={1\over 2}\onesub y{ik}y^{kl}\onesub y{lj}
-2\dddd {\b R}0i0j$$
$$y_{ij}(0)=\rho_{ij}~~~;~~~\onesub y{ij}(0)=2\lambda_{ij}~.$$
{}From the theory of ODE it follows therefore that $~\sigma_{ij}=a_{ij}~$.
For (b) we first look at the tensor $~\rho^{kl}(\onesub {\eta}{ki}\eta
_{lj}-\eta_{ki}\onesub {\eta}{lj})~$ and by differentiating this we see that
it vanishes identically. Thus $~\rho^{kl}\onesub {\eta}{ki}\eta_{lj}~$ is
symmetric in $~i~$ and $~j~$, then (b) is obtained by differentiation of
$~\sigma_{ij}~$. (c) is a direct consequence of (a), since locally
$$~d\mu_t=
\sqrt{\det(\sigma_{ij})}dx=\sqrt{\det(\rho^{ik}\eta_{kj})^2\det(\rho_{ij})}dx~
=\det(\rho^{ik}\eta_{kj})d\mu_0~,$$
where the last equality is true since $~\eta_{ij}~$ is invertible and
therefore positive definite ($~\eta_{ij}(0)~$ is positive definite).
Furthermore
$~\det(\rho^{ik}\eta_{kj})~$ does not depend on
the choice of a coordinate system.
\bigskip\noindent
For a constant $~\kappa~$ define the function
$$s_\kappa:=\cases{{1\over\sqrt{\kappa}}\sin(\sqrt{\kappa}r)&$\kappa>0~;$\cr
\cr
r&$\kappa=0~;$\cr
\cr
{1\over\sqrt{-\kappa}}\sinh(\sqrt{-\kappa}r)&$\kappa<0~.$\cr}$$
Let $~sym_k~$ be the $~k$-th elementary symmetric function, i.e.
for all $~a_1,\dots a_n\in\real {}~$
$$sym_k(a_1,\dots,a_n):=
\sum_{1\le i_1<\cdots <i_k\le n}a_{i_1}\cdots a_{i_k},~~~~0\le k\le n$$
and set
$$S_k:=sym_k(\lambda_1,\dots,\lambda_n)~,$$
where $~\lambda_1,\dots,\lambda_n~$ are the principal curvatures on $~M_0~$.

{\corollar {1.3~(Steiner~formula)} \itemitem{}
If the sectional curvature $~\kappa~$ is constant
and if we denote the area of $~M_r~$ by $~\vert M_r\vert~$, then
$$\eta_{ij}=\onesub s{\kappa}\rho_{ij}+s_\kappa\lambda_{ij}$$
$$\vert M_r\vert =\sum_{l=0}^n\onesub s{\kappa}^{l}s_\kappa^{~n-l}\int_{M^n}
S_{\scriptscriptstyle {n-l}}(0)d\mu_0~.$$
}
\proof That $~\eta_{ij}~$ satisfies the above equation follows from the fact
that it satisfies the appropriate evolution equation. The area formula is
then a consequence of Theorem 1.2 (c).
\remark
For an excellent book on Weyl's tube formula and Steiner type formulas see [G].

\bigskip\noindent
{\prop {1.4}\itemitem{(a)} Assume that $~\dddd {\b R}0i0j
\le\mu(r) \sigma_{ij}~$,
$~\lambda_{ij}\ge b_l\rho_{ij}~$ and that $~f_\mu~$ is any
function with $~\twosub f{\mu}+\mu f_\mu\le0~$, $~f_\mu(0)=1~$
and $~\onesub f{\mu}(0)\le b_l~$. Then
$$\tau_{ij}\ge{\onesub f{\mu}\over f_\mu}\sigma_{ij}~,$$
as long as $~f_\mu>0~$.
\smallskip\noindent\itemitem{(b)} Assume that $~\dddd {\b R}0i0j
\ge\nu(r) \sigma_{ij}~$,
$~\lambda_{ij}\le b_u\rho_{ij}~$ and that $~f_\nu~$ is any
function with $~\twosub f{\nu}+\nu f_\nu\ge0~$, $~f_\nu(0)=1~$
and $~\onesub f{\nu}(0)\ge b_u~$. Then
$$\tau_{ij}\le{\onesub f{\nu}\over f_\nu}\sigma_{ij}~,$$
as long as $~f_\nu>0~$.
$$$$
}
\proof
We prove only part (a) since the proofs for (a) and (b) are almost the same.
Assume that for positive $~\epsilon$,
$~f_\epsilon~$ is the solution of the Jacobi equation
$$\twosub f{\epsilon}
=-(\mu+\epsilon) f_\epsilon~~~,~~~f_\epsilon(0)=1~~~,~~~\onesub f{\epsilon}
(0)=b_l-\epsilon~,$$
and let $~c_\epsilon
:={\onesub f{\epsilon}\over f_\epsilon}~$. Define the symmetric tensor
$~M_{ij}:=\tau_{ij}-c_\epsilon \sigma_{ij}~$. Assume that $~r_0~$ is
the first distance
where a zero eigenvector $~v~$ of $~M_{ij}~$ at one point $~x\in M~$
occurs and let us assume that at this point we have chosen normal coordinates
with respect to $~\sigma_{ij}~$
such that $~M_{ij}~$ (and therefore also $~\tau_{ij}~$)
becomes diagonal and we can also assume (after a rotation of the coordinate
frame) that at this point $~v^i=\delta ^i_1~$. Now since $~r_0~$ is the first
distance where $~M_{ij}v^iv^j~$ vanishes we must have
$~\td {}r{}_{\vert r_0}M_{ij}v^iv^j\le 0~$.
But on the other hand we calculate from Lemma 1.1 that
$$\td {}r{}_{\vert r_0} M_{ij}v^iv^j=(\tau_{il}\sigma^{lm}\tau_{mj}-\dddd
{\b R}0i0j-2c_\epsilon
\tau_{ij}-\onesub c{\epsilon} \sigma_{ij})v^iv^j$$
$$=\tau_{11}^2-\dddd {\b R}0101-2c_\epsilon \tau_{11}-\onesub c{\epsilon}$$
$$=-(\onesub c{\epsilon}+c_\epsilon^2+\dddd {\b R}0101)>-(\onesub c{\epsilon}
+c_\epsilon^2+\mu+\epsilon)=0~.$$
This contradiction proves that $~\tau_{ij}>{\onesub f{\epsilon}
\over f_\epsilon}\sigma_{ij}~$
as long as $~f_\epsilon >0~$. Since the solution of the above Jacobi equation
depends continously on $~\epsilon~$ we have
$$\tau_{ij}\ge{\onesub f0\over f_0}\sigma_{ij}~,$$
as long as $~f_0>0~$, where
$$\twosub f0=-\mu f_0~~~,~~~f_0(0)=1~~~,~~~\onesub f0(0)=b_l~.$$
Now let $~\tau_\mu:=\max\{r\in(0,R)\vert f_\mu(r)>0\}~$,
$\tau_0:=\max\{r\in(0,R)\vert f_0(r)>0\}~$.
Define the function $~h:=f_\mu\onesub f0-\onesub f{\mu} f_0~$. Then we have
$~h(0)\ge 0~$ and $~\one h=\twosub f0 f_\mu-f_0\twosub f{\mu}
\ge 0~$ for $~r\in\min\{\tau_0,\tau_\mu\}~$. Thus $~\left ({f_0\over
f_\mu}\right)\acute {}={h\over f_\mu^2}\ge 0~$ on
$~[0, \min\{\tau_0,\tau_\mu\})~$.
But since $~{f_0\over f_\mu}(0)=1~$, this implies that $~f_0\ge f_\mu~$ on
$~[0,\min\{\tau_0,\tau_\mu\})~$ and thus $~\tau_0\ge\tau_\mu~$. So
$~h(r)\ge 0~$ for $~r\in[0, \tau_\mu)~$ and consequently also
${\onesub f0\over f_0}\ge{\onesub f{\mu}\over f_\mu}~$ for
$~r\in[0,\tau_\mu)~$.
This proves the proposition.

{\corollar {1.5}\itemitem{}
$$f_\nu^2\rho_{ij}\ge \sigma_{ij}$$
as long as $~f_\nu>0~$ and
$$\sigma_{ij}\ge f_\mu^2 \rho_{ij}$$
as long as $~f_\mu>0~$.
\smallskip\noindent
}
\proof
We have $~\td {}r{}({1\over f_\mu^2}\sigma_{ij}-\rho_{ij})
={2\over f_\mu^2}(\tau_{ij}-
{\onesub f{\mu}\over f_\mu}\sigma_{ij})\ge 0~~$ and $~~\td {}r{}({1\over
f_\nu^2}\sigma_{ij}-\rho_{ij})\le 0~$.

{\corollar {1.6~(Rauch~comparison~theorem)}\itemitem{}
Let $~N_1,~N_2~$ be two Riemannian manifolds with sectional curvatures
$~\kappa_1,~\kappa_2~$ and let $~c_1,c_2~$ be two geodesics parametrized
by arclength and assume that along these geodesics the inequality
$~\kappa_1\ge \kappa_2~$ holds. Assume $~J_1,~J_2~$ are two Jacobi fields
along these geodesics such that $~\langle J_i(0),\onesub ci(0)\rangle=
\langle \onesub Ji(0),\onesub ci(0)\rangle=0~$
and $~\vert J_1(0)\vert =
\vert J_2(0)\vert,~\vert \onesub J1(0)\vert=\vert \onesub J2
(0)\vert~$. Then as
long as there are no conjugate points on $~c_1~$ we have the
estimate
$$\vert J_1(r)\vert\le\vert J_2(r)\vert~.$$
\smallskip\noindent}
\proof
First assume that $~w_i:=J_i(0)\neq 0~$. Let $~M_i~$ be a small hypersurface,
normal to the geodesic $~c_i~$ such that $~\b \nabla^{\top}
_{w_i}\nu_i=v_i:=\onesub Ji(0)~$
(it is clear that such surfaces exist). Furthermore let $~\gamma_i(s)~$ be
curves on $~M_i~$ such that $~\gamma_i(0)=c_i(0),~\onesub {\gamma}i(0)=w_i~$
and define
a family of geodesics $~c_i(s,r):=\ex {\gamma_i(s)}{r\nu_i(s)}$,
where $~\nu_i(s)~$
is that
unit normal vector field along $~\gamma_i$ with $~\nu_i(0)=\onesub ci(0)~$.
Then we have $~J_i(r)={\partial\over\partial s}_{\vert s=0}c_i(s,r)~$.
Now choose $~\kappa~$ such that $~\kappa_1\ge\kappa\ge\kappa_2~$ and solve the
Fermi flow for both $~M_i~$. Using Corollary 1.5 we see that
$~\vert J_1(r)\vert^2\le f_\kappa^2\vert J_1(0)\vert^2=f_\kappa^2\vert J_2(0)
\vert^2\le\vert J_2(r)\vert^2~$ as long as both solutions exist. That shows
that no conjugate point on $~c_2~$ can occur before one occurs on $~c_1~$. If
we now choose the initial hypersurfaces smaller and smaller, we see that the
maximal time for which both solution exist tends to the maximal time
for which there are no conjugate points on $~c_1~$.
\smallskip\item{}
In the case where $~w_i=0~$ we set
$~\t w_i:=\epsilon z_i~$ for some arbitrary fixed  unit vectors $~z_i~$ and let
$~\epsilon ~$ tend to zero. This gives the result since the solution of
the Jacobi equation depends continously on initial data.

\bigskip\noindent
The evolution equations for the mean curvature and the volume form
immediately imply an estimate for the volume of a geodesic ball, if
the Ricci curvature is bounded below. This well known result is due to
Heintze and Karcher (see for example [GHL]). Unfortunately, almost no results
are known when the Ricci curvature is bounded above. We are able
to prove such a result under an additional assumption on the curvature tensor.

{\defin {1.7}\itemitem{}
We say that the sectional curvature $~\b\sigma~$
is $~\epsilon-$Ricci pinched, if for any vector $~X~$ and any two-plane
$~X\wedge Y~$ the estimate $~\vert
\b\sigma(X\wedge Y)-{1\over n}\b Ric(X,X)\vert\le \epsilon~$ holds, where
the dimension of $~N~$ is $~n+1~$.
\smallskip\noindent}
\bigskip\noindent
We are now going to show, that the $~\epsilon-$Ricci pinching is strongly
related to the umbilicness of a convex surface.

{\prop {1.8}\itemitem{}
Assume $~\b\sigma~$ is $~\epsilon-$Ricci pinched and $~M_0~$ is the immersion
of a convex hypersurface such that $~\vert \tau_{ij}-{\tau\over n}\sigma_{ij}
\vert\le c\sigma_{ij}~$.
Then as long as $~M_r~$ stays convex, we have the estimate
$$\vert \tau_{ij}-{\tau\over n}\sigma_{ij}\vert \le (c+\epsilon r)\sigma_{ij}~.
$$}
\proof
Define the tensor
$$M_{ij}:=\tau_{ij}-{\tau\over n}\sigma_{ij}+(c+\t\epsilon r)\sigma_{ij}$$
and
$$N_{ij}:=\tau_{ij}-{\tau\over n}\sigma_{ij}-(c+\t\epsilon r)\sigma_{ij}~,$$
where $~\t\epsilon>\epsilon~$
As in the proof of Proposition 1.4 we look at the first time where a zero
eigenvector of $~M_{ij}~$ or $~N_{ij}~$ occurs. With the conventions
as above assume that $~v~$ is a zero eigenvector of $~N_{ij}~$ (the case
$~M_{ij}~$ will be analogous). Then we have
$$\td {}r{}_{\vert r_0} N_{ij}v^iv^j=(\tau_{il}\sigma^{lm}\tau_{mj}-\dddd
{\b R}0i0j-2({\tau\over n}+c+\t\epsilon r)\tau_{ij}+{1\over n}(\vert
\tau_{ij}\vert^2+\b Ric(\nu,\nu))\sigma_{ij}-\t\epsilon\sigma_{ij})\u vi\u vj$$
$$=-\tau_{11}^2+{1\over n}\vert \tau_{ij}\vert^2-\dddd {\b R}0101+{1\over n}
\b Ric(\nu,\nu)-\t\epsilon$$
If $~M_r~$ is still convex, we have $~\vert\tau_{ij}\vert^2\le n\tau_{11}^2~$
and therefore using the $~\epsilon-$Ricci umbilicness
$$\td {}r{}_{\vert r_0} N_{ij}v^iv^j\le \epsilon-\t\epsilon<0$$
which is a contradiction. For $~\t\epsilon\to\epsilon~$ we get the desired
result.

\bigskip\noindent
{\prop {1.9}\itemitem{}
Let $~N~$ be $~\epsilon-$Ricci pinched with $~\b Ric\le d~$. Let
$~m:={d\over n}-\epsilon~$ and define $~R_0:={\pi\over 2\sqrt{{d\over
n}+\epsilon}}~$, if $~{d\over n}+\epsilon>0~$
and otherwise set $~R_0:=\infty~$
{}.
Then any geodesic sphere of radius smaller than $~R_0~$ is convex,
and if we assume that $~S_{r_1},~S_{r_2}~$ are two
geodesic spheres with $~0<r_1<r_2<R_0~$ such that on $~S_{r_1}~$
$~\tau>\tau_-~$ and $~\vert \tau_{ij}-{\tau\over n}\sigma_{ij}\vert\le c\sigma
_{ij}~$,
then we can estimate the enclosed volume between these spheres,
as long as the function $~u(r):=({\tau_-\over n}+c)s_m(r)+\onesub sm(r)~$
is positive, by
$$V(r_1,r_2)\ge\vert S_{r_1}\vert\int_0^{r_2-r_1}u^n(r)e^{-ncr-
{n\epsilon\over 2}r^2}dr~.$$
\smallskip\noindent}
\proof
The assumptions on the curvature imply that the sectional curvature is bounded
between $~{d\over n}-\epsilon<\kappa<{d\over n}+\epsilon~$. Using
Proposition 1.4 with $~\mu={d\over n}+\epsilon~$
one easily checks that any geodesic sphere
of radius smaller than $~R_0~$ is convex since all principal curvatures
tend to infinity as the radius tends to zero, and therefore $~b_l~$ can be
chosen to be arbitrarily large. Now we use Proposition 1.8 and obtain
$$\vert \tau_{ij}\vert^2\le n({\tau\over n}+c+\epsilon r)^2~.$$
Define $~y:={\tau\over n}+c+\epsilon r~$. Then we calculate
$$\one y=-{1\over n}(\vert \tau_{ij}\vert^2+\b Ric(\nu,\nu))+\epsilon
\ge -y^2-m~.$$
Let $~f~$ be the solution of $~\one f=-f^2-m~,f(0)={\tau_-\over n}+c~$, i.e.
$~f=\one {(\ln u)}~$ with $~u~$ as above. Then $~\one {(y-f)}\ge -(y+f)(y-f)~$
and consequently $~y\ge f~$ as long as $~r_1+r <R_0~$ and $~f~$
is defined, i.e. $~u~$ is positive. This implies
$$\one {\vert S_{r_1+r}\vert}=\int \tau d\mu\ge n\vert S_{r_1+r}
\vert\one {(\ln u-cr-{\epsilon\over 2}r^2)}~,$$
thus
$$\vert S_{r_1+r}\vert\ge \vert S_{r_1}
\vert u^n(r)e^{-ncr-n{\epsilon\over 2}r^2}~.$$
Integration from $~0~$ to $~r_2-r_1~$ gives the result.


\section {2. Short- and longtime existence results for the Fermi flow}
\smallskip\noindent
So far we have not proven that a solution exists for some short
time interval. In this section we prove the existence of a unique solution
and also derive estimates for the maximal time interval on which
it is admitted. In particular we show that in some cases
we obtain eternal solutions.

{\prop {2.1}\itemitem{}
Let $~F_0:M^n\to N^{n+1}~$ be a smooth immersion of an orientable and compact
manifold and let $~\nu~$ be one of the unit normal vector fields on $~M~$.
Then there exists an $\epsilon >0$,
a unique family of smooth submanifolds $~M_r~$
and a smooth family of local diffeomorphisms $~F_r:M^n\to N^{n+1}~$
such that
$$\td {F_r}r{}=\nu~~~,~~~M_r=F_r(M^n)~~~\forall~ r\in[0,\epsilon)$$
}
\proof
Let $~c_x(r):=\ex {F_0(x)}{r\nu(F_0(x))}~$ and define $~F_r(x):=c_x(r)~$. This
is well defined since $~M^n~$ is compact and therefore the geodesic equation
can be solved uniformly for all $~x\in M^n~$ on some short time interval
$~[0,\epsilon)~$. Moreover, $~F_r~$ is smooth since the solution of the
geodesic equation depends smoothly on the initial conditions.
If we assume $~\epsilon~$
to be so small that the exponential map at each point on the initial
hypersurface is injective, then $~F_r~$ becomes a local diffeomorphism.
{}From the
Gauss lemma it follows that $~\td {F_r}r{}=\td {c_x}r{}~$ is a unit normal to
the hypersurfaces defined by $~F_r~$. This proves existence. Uniqueness
follows from Lemma 1.1(c) and
the fact that the solution of the geodesic equation is unique.

\remark
We have just proved that it is possible to reduce the nonlinear system of
first order PDE's that is defined through (1) to an n-parameter family
of nonlinear systems of second order ODE's modeled over $~M^n~$.

{\defin {2.2}\itemitem{}
Let $~M~$ be an embedded, orientable hypersurface with unit normal vector
field $~\nu~$ and let $~\Gamma~$ be the set of all geodesics $~\gamma~$
parametrized by arclength with $~\gamma(0)\in M,~\one\gamma(0)=\nu~$.
The self distance $~d_{\nu}~$
of $~M~$ with respect to $~\nu~$ is defined as
$$d_{\nu}:=\sup\{\rho\vert \Gamma_{\vert (0,\rho)}\cap M=\emptyset\}~.$$
\smallskip\noindent}

{\lemma {2.3}\itemitem{}
Assume that $~M_0~$ is the embedding of a closed, orientable hypersurface
in a Riemannian manifold $~N~$ such that $~N~$ is connected and $~N-M_0~$ is
disconnected and consists of two components $~A~$ and $~B~$ and let $~d_{\nu}~$
be the self distance of $~M_0~$ with respect to the normal vector field that
points in direction of $~A~$. Then $~{d_{\nu}\over 2}~$ is the first time
where $~M_t~$ fails to be an embedding (provided the flow exists this far)
and for all
$~t\le{d_{\nu}\over 2}~$ and all $~p\in M_t~$ we have $~d(p,M_0)=t~$.
\smallskip\noindent}
\proof
Let us first remark that since $~M_0~$ is compact and smooth,
it is clear that there exists an $~\epsilon>0~$
such that $~d(p,M_0)=t~$ for all $~t\le \epsilon~$ and all $~p\in M_t~$.
\smallskip\item{}
Now let $~T~$ be the first time where this property fails to be true,
i.e. the first time such that there exists at least one point
$~p\in M_T~$ and two distance-minimizing geodesics $~\gamma_1,~\gamma_2~$
between $~M_T~$ and $~M_0~$, parametrized by arclength with $~\gamma_i(0)=p~$
and $~\gamma_i(T)\in M_0~$. These geodesics are both normal to both surfaces
and moreover $~\onesub {\gamma}i(T)=-\nu(\gamma_i(T))~$.
This follows from the above remark
and the fact the $~N-M_0~$ is disconnected; therefore,
$$d(p,M_0)=\inf\{l(\gamma)\vert\gamma\subset N {\rm~connects~}
p{\rm~with~} M_0\}$$
$$=\inf \{l(\gamma)\vert \gamma \subset A
{\rm~connects~}p{\rm~with~} M_0\}~.$$
Thus $~\gamma_1\cup\gamma_2~$ is a smooth geodesic arc connecting $~M_0~$
with itself. Therefore $~d_{\nu}\le 2T~$. Since both geodesics $\gamma_i~$
satisfy $~\onesub {\gamma}i(T)=-\nu(\gamma_i(T)),~\gamma_i(0)=p~$, this also
shows that $~M_T~$ touches itself and cannot be more than immersed.
On the other hand if $~\t T<T~$ were the first time where
$~M_{\t T}~$ becomes immersed, then $~M_{\t T}~$ would touch itself at
some point $~p\in M_{\t T}~$, and from the definition
of $~T~$ we would obtain a contradiction.

{\prop {2.4}\itemitem{}
Let $~N^{n+1}~$ be geodesically complete and connected and let
$~F_0:M^n\to N^{n+1}~$ be a smooth embedding of an orientable closed
manifold such that $~N^{n+1}-M^n~$ is disconnected and consists of two
components $~A~$ and $~B~$. Assume that $~\nu~$ is a unit normal vector
field on $~M_0~$ that points in direction of $~A~$, that $~d_{\nu}~$ is the
self distance of $~M_0~$, and that at
any point $~y\in ~A~$ the sectional curvature $~\sigma(y)~$ is
bounded above or below by functions
$\mu(d(y))={c_{\mu}\over (1+a_{\mu}d(y))^2}~$, or
$\nu(d(y))={c_{\nu}\over (1+a_{\nu}d(y))^2}~$ respectively,
where $~a_{\mu},~a_{\nu}>0~$,
$~c_{\mu},~c_{\nu}~$ are
fixed constants and $~d(y)~$ is the distance to $~M_0~$.
Further, let $~b_l,~b_u~$ be bounds for the second fundamental form on
$~M_0~$, i.e. $~b_u\rho_{ij}\ge\lambda_{ij}\ge b_l\rho_{ij}~$,
and let $~m_i:={4c_i\over a_i^2}-1~$ for $~i=\mu, \nu~$
Then we have, with $~\vert a,b\vert:=\min\{a,b\}~$,
\smallskip\itemitem {(a)} $~m_\mu\le 0~$ and
$~{2\over a_\mu}b_l-1\ge -\sqrt{-m_\mu}\Rightarrow r_0=\vert
{d_{\nu}\over 2},\infty\vert~$
\smallskip\itemitem {(b)} $~m_\mu=0~$ and
$~{2\over a_\mu}b_l-1< -\sqrt{-m_\mu}\Rightarrow
r_0\ge \vert {d_{\nu}\over 2},{1\over a_\mu}(e^{{2a_\mu\over
a_\mu-2b_l}}-1)\vert~$
\smallskip\itemitem {(c)} $~m_\mu< 0~$ and
$~{2\over a_\mu}b_l-1< -\sqrt{-m_\mu}\Rightarrow
r_0\ge\vert {d_{\nu}\over 2}
, {1\over a_\mu}(\root {\sqrt{-m_\mu}} \of
{ {a_\mu-2b_l+a_\mu\sqrt{-m_\mu}\over a_\mu
-2b_l-a_\mu\sqrt{-m_\mu}} }-1)\vert~$
\smallskip\itemitem {(d)} $~m_\mu>0~$ and
$~{2\over a_\mu}b_l-1=0\Rightarrow r_0\ge\vert {d_{\nu}\over 2}
,{1\over a_\mu}(e^{\pi\over\sqrt{m_\mu}}-1)\vert$
\smallskip\itemitem {(e)} $~m_\mu>0~$ and
$~{2\over a_\mu}b_l-1<0\Rightarrow r_0\ge\vert {d_{\nu}\over 2}
,{1\over a_\mu}(e^{{2\over\sqrt{m_\mu}}
\arctan {a_\mu\sqrt{m_\mu}\over a_\mu-2b_l}}-1)\vert$
\smallskip\itemitem {(f)} $~m_\mu>0~$ and
$~{2\over a_\mu}b_l-1>0\Rightarrow r_0\ge\vert {d_{\nu}\over 2}
,{1\over a_\mu}(e^{{2\over\sqrt{m_\mu}}
(\pi+\arctan {a_\mu\sqrt{m_\mu}\over a_\mu-2b_l})}-1)\vert$
\smallskip\itemitem {(g)} $~m_\nu=0~$ and
$~{2\over a_\nu}b_u-1< -\sqrt{-m_\nu}\Rightarrow
r_0\le {1\over a_\nu}(e^{{2a_\nu\over a_\nu-2b_u}}-1)~$
\smallskip\itemitem {(h)} $~m_\nu< 0~$ and
$~{2\over a_\nu}b_u-1< -\sqrt{-m_\nu}\Rightarrow
r_0\le {1\over a_\nu}(\root {\sqrt{-m_\nu}} \of
{ {a_\nu-2b_u+a_\nu\sqrt{-m_\nu}\over a_\nu
-2b_u-a_\nu\sqrt{-m_\nu}} }-1)~$
\smallskip\itemitem {(i)} $~m_\nu>0~$ and
$~{2\over a_\nu}b_u-1=0\Rightarrow
r_0\le{1\over a_\nu}(e^{\pi\over\sqrt{m_\nu}}-1)$
\smallskip\itemitem {(j)} $~m_\nu>0~$ and
$~{2\over a_\nu}b_u-1<0\Rightarrow r_0\le{1\over a_\nu}(e^{{2\over\sqrt{m_\nu}}
\arctan {a_\nu\sqrt{m_\nu}\over a_\nu-2b_u}}-1)$
\smallskip\itemitem {(k)} $~m_\nu>0~$ and
$~{2\over a_\nu}b_u-1>0\Rightarrow r_0\le{1\over
a_\nu}(e^{{2\over\sqrt{m_\nu}}(
\pi+\arctan {a_\nu\sqrt{m_\nu}\over a_\nu-2b_u})}-1)~.$
}
\proof
First we obtain from lemma 2.3 that $~d(F_r(x))=r,~\forall ~x\in M^n~$, as long
as $~r\le {d_{\nu}\over 2}~$.
Now we can use Proposition 1.4 with $~\mu(r)={c_\mu\over (1+a_\mu r)^2}~$,
or $~\nu(r)={c_\nu\over (1+a_\nu r)^2}~$ resp. Since $~N~$ is geodesically
complete, the maps $~F_r(x):=\ex {F_0(x)}{r\nu(F_r(x))}~$ are welldefined, and
as in the proof of proposition 2.1 it follows that this is a smooth family
of maps and a family of diffeomorphisms provided $~r~$ is small enough. So the
only thing that can go wrong is that $~\sigma_{ij}~$ could degenerate, i.e.
the eigenvalues of $~\sigma_{ij}~$ could tend to infinity or zero.
Using Proposition 1.4 and the fact that the solution of the Jacobi equation
$~\two f+{c\over (1+ar)^2}f=0~,f(0)=1,~\one f(0)=b$ is explicitly given by
$$f(r)=\sqrt{1+ar}\bigl (({2\over a}b-1)s_m\ln \sqrt{1+ar}+\onesub sm\ln\sqrt
{1+ar}\bigr)~,$$
with $~m:={4c\over a^2}-1~$, we see that no eigenvalue of the metric can tend
to infinity in finite time. So an estimate for the distance, where $~f_\mu,~
f_\nu~$ become zero for the first time, gives lower and upper bounds for
$~r_0~$. One easily checks that these estimates are given by those stated in
the proposition.

{\prop {2.5} \itemitem{}
Under the assumption of Proposition 2.4 assume that a lower bound
for the sectional curvature is now given by $~\nu(r)={c\over (1+ar)^{2-
\epsilon}}~$
, where $~a,~c,~\epsilon >0~$. Then there exists no eternal solution, i.e.
$~r_0<\infty~$.
\smallskip\noindent}
\proof
Assume that $~r_0=\infty~$ and choose $~R_0~$ so that
$~4c(1+aR_0)^\epsilon>a^2~$.
Then we define $~\t c:=c(1+aR_0)^\epsilon~$ and conclude that for all
$~r\ge R_0~$
$$\nu(r)={c\over (1+ar)^2}(1+ar)^\epsilon\ge {\t c\over (1+ar)^2}~.$$
Using Proposition 2.4 we see that since $~\t m_\nu>0~$,
a further extension of the flow
can only be made for some finite time, contradicting to $~r_0=\infty~$.

\remark
Proposition 2.5 means that a manifold cannot admit an asymptotically
flat end if the curvature does not decay fast enough. One could also try to
give estimates for $~r_0~$ in the case, where the functions $~\mu,~\nu~$
are given by $~{c\over (1+ar)^\alpha}~$ with some constant $~\alpha~$.
However, it turns out that the corresponding Jacobi equation admits a much
more complicated solution (see [K], $2\cdot 14$). By Proposition 1.4 it would
be enough to give estimates for sub- or supersolutions of the Jacobi equation.

\bigskip\noindent
Let us finally prove a statement concerning the gradient growth of the second
fundamental form

{\prop {2.6}\itemitem{}Assume that $~{\onesub f{\mu}\over f_\mu}>0~$ on [0,R)
and that there exist positive constants $~c_1,~c_2~$ such that
$~\vert\b\nabla_k\dddd {\b R}0i0j\vert^2\le c_1\bigl(
{\onesub f{\mu}\over f_\mu}\bigr) ^2~$
and
$~\vert\langle\tau_{kp},\dddd {\b R}0jpi\rangle\vert^2\le c_2
\bigl({\onesub f{\mu}\over f_\mu}\bigr)^2~$.
Then for any $~\epsilon>0~$ there exist constants
$~a_1,~a_2~$ such that for all $~r\in[0,R)~$
$$\vert\cd {\dd {\tau}ij}k\vert^2\le a_1f_{\mu}^{5\epsilon-6}+a_2~,$$
with $~a_1~$ and $~a_2~$ given by
$$a_2={c_1+4c_2\over\epsilon(6-5\epsilon)},~a_1=\min_{r=0}
\vert\cd {\dd {\tau}ij}k\vert^2-a_2~.$$
}
\proof
First we calculate
$$\td {}r{} \nabla_k\tau_{ij}=\nabla_k(\tau_{iu}\sigma^{uv}\tau_{vj}-
\dddd {\b R}0i0j)-\tau_{pj}\sigma^{pq}(\nabla_k\tau_{qi}+\nabla_i\tau_{qk}
-\nabla_q\tau_{ki})$$
$$-\tau_{pi}\sigma^{pq}(\nabla_k\tau_{qj}+\nabla_j\tau_{qk}
-\nabla_q\tau_{kj})$$
$$=-\cdb {\dddd {\b R}0i0j}k-\tau_{kp}\sigma^{pq}\dddd {\b R}qi0j
-\tau_{kp}\sigma^{pq}\dddd {\b R}0iqj$$
$$-\sigma^{pq}(\tau_{pj}\nabla_i\tau_{qk}-\tau_{pj}\nabla_q\tau_{ki}
+\tau_{pi}\nabla_j\tau_{qk}-\tau_{pi}\nabla_q\tau_{kj})$$
$$=-\cdb {\dddd {\b R}0i0j}k-\sigma^{pq}(\tau_{kp}\dddd {\b R}0jqi
+\tau_{kp}\dddd {\b R}0iqj+\tau_{pj}\dddd {\b R}0kqi
+\tau_{pi}\dddd {\b R}0kqj)$$
where we have used the Codazzi equations
$~\nabla_i\tau_{kq}-\nabla_q\tau_{ki}=
\dddd {\b R}0kqi~$
(Note that $~\td {}r{} \nabla_k\tau_{ij}=0~$
in the case where $~N~$ is a space form), and then
$$\td {}r{} \vert\cd {\dd {\tau}ij}k\vert^2=-2\uu {\sigma}ks\uu {\sigma}lt
\uu {\sigma}in\uu {\sigma}jm\tau_{st}\nabla_k\tau_{ij}\nabla_l\tau_{nm}
-4\uu {\sigma}kl\uu {\sigma}is
\uu {\sigma}nt\uu {\sigma}jm\tau_{st}\nabla_k\tau_{ij}\nabla_l\tau_{nm}$$
$$+2\sigma^{kl}\sigma^{in}\sigma^{jm}\nabla_l\tau_{mn}\td
{}r{}\nabla_k\tau_{ij}$$
$$=-2\uu {\sigma}ks\uu {\sigma}lt
\uu {\sigma}in\uu {\sigma}jm\tau_{st}\nabla_k\tau_{ij}\nabla_l\tau_{nm}
-4\uu {\sigma}kl\uu {\sigma}is
\uu {\sigma}nt\uu {\sigma}jm\tau_{st}\nabla_k\tau_{ij}\nabla_l\tau_{nm}$$
$$-2(\cdb {\dddd {\b R}0i0j}k+2\tau_{kp}\sigma^{pq}\dddd {\b R}0jqi+2\tau_{pj}
\sigma^{pq}\dddd {\b R}0kqi)\sigma^{kl}\sigma^{in}
\sigma^{jm}\nabla_l\tau_{mn}~.$$
Using Schwartz's inequality we obtain, for any positive $~\eta~$,
$$\td {}r{}\vert\cd {\dd {\tau}ij}k\vert^2\le -2\uu {\sigma}ks\uu {\sigma}lt
\uu {\sigma}in\uu {\sigma}jm\tau_{st}\nabla_k\tau_{ij}\nabla_l\tau_{nm}
-4\uu {\sigma}kl\uu {\sigma}is
\uu {\sigma}nt\uu {\sigma}jm\tau_{st}\nabla_k\tau_{ij}\nabla_l\tau_{nm}$$
$$+\eta^{-1}(\vert\b\nabla_k\dddd {\b R}0i0j\vert^2+4\vert\langle\tau_{kp},
\dddd {\b R}0jpi\rangle\vert^2)+5\eta\vert\nabla_k\tau_{ij}\vert^2~.$$
Now we use Proposition 1.4 (a) and obtain
$$\td {}r{}\vert\cd {\dd {\tau}ij}k\vert^2\le -6{\onesub f{\mu}\over f_\mu}
\vert\cd {\dd {\tau}ij}k\vert^2
+\eta^{-1}(\vert\b\nabla_k\dddd {\b R}0i0j\vert^2+4\vert\langle\tau_{kp},
\dddd {\b R}0jpi\rangle\vert^2)+5\eta\vert\nabla_k\tau_{ij}\vert^2~,$$
and if we choose $~\eta=\epsilon{\onesub f{\mu}\over f_\mu}~$, then we get
$$\td {}r{}\vert\cd {\dd {\tau}ij}k\vert^2\le
(5\epsilon-6){\onesub f{\mu}\over f_\mu}
\vert\cd {\dd {\tau}ij}k\vert^2+{c_1+4c_2\over\epsilon}
{\onesub f{\mu}\over f_\mu}~.$$
Thus we have shown
$$\td {}r{}\Bigl(\ln\bigl((\vert\cd {\dd
{\tau}ij}k\vert^2-a_2)f_{\mu}^{6-5\epsilon}\bigr)\Bigr)\le 0~,$$
and consequently the result.

\bigskip\noindent
The same calculations prove the following
{\prop {2.7}\itemitem{}Assume that $~{\onesub f{\mu}\over f_\mu}<0~$ on [0,R)
and that there exist positive constants $~c_1,~c_2~$ such that
$~\vert\b\nabla_k\dddd {\b R}0i0j\vert^2\le c_1\bigl(
{\onesub f{\mu}\over f_\mu}\bigr) ^2~$
and
$~\vert\langle\tau_{kp},\dddd {\b R}0jpi\rangle\vert^2\le c_2
\bigl({\onesub f{\mu}\over f_\mu}\bigr)^2~$.
Then for any $~\epsilon<0~$ there exist constants
$~a_1,~a_2~$ such that for all $~r\in[0,R)~$
$$\vert\cd {\dd {\tau}ij}k\vert^2\le a_1f_{\mu}^{5\epsilon-6}+a_2~,$$
with $~a_1~$ and $~a_2~$ given by
$$a_2={c_1+4c_2\over\epsilon(6-5\epsilon)},~a_1=\min_{r=0}
\vert\cd {\dd {\tau}ij}k\vert^2-a_2~.$$.
}


\section {References}
\bigskip\bigskip\noindent
\cite {G}{Gray, A.;}{Tubes,}{Addison Wesley,}{1990}{}
\cite {GHL}{Gallot, S.; Hulin, D.; Lafontaine, J;}{Riemannian Geometry,}
{Springer Verlag,}{Berlin, New York,}{1987}
\cite {K}{Kamke;}{Differentialgleichungen, L\"osungsmethoden und L\"osungen,}
{Chelsea Publishing Company,}{New York}{}
\cite {S}{Smoczyk, K.;}{Symmetric hypersurfaces in Riemannian manifolds
contracting to Lie groups by their mean curvature,}{Calc. Var.}{1996, 4 (02),}
{P. 155}

\end